\documentclass[11pt,a4paper]{article}
\usepackage{a4wide}
\usepackage{latexsym}
\usepackage{amsmath}
\usepackage{amsthm}
\usepackage{amssymb,enumerate}
\usepackage[usenames]{color}
\usepackage{graphicx}
\usepackage{fancybox}
\usepackage{hyperref}
\usepackage[latin1]{inputenc}
\theoremstyle{plain}
\newtheorem{theo}{Theorem}[section]
\newtheorem{lemma}[theo]{Lemma}
\newtheorem{prop}[theo]{Proposition}

\theoremstyle{definition}

\newenvironment{proof1}{\medskip\par\noindent{\bf Proof.}}{\hfill $\Box$
\medskip\par}

\def\C{\mathbb{C}}
\def\N{\mathbb{N}}
\def\R{\mathbb{R}}

\def\a{\alpha}
\def\b{\beta}
\def\o{\omega}

\def\ga{\gamma}
\def\Ga{\Gamma}

\def\bM{\mathbb{M}}
\def\M{\mathbb{M}}
\def\hM{\widehat{\M}}

\def\bL{\mathbb{L}}
\def\bm{\boldsymbol{m}}
\def\m{\boldsymbol{m}}
\def\bl{\boldsymbol{\ell}}

\definecolor{azulosc}{rgb}{0.2,0.1,0.7}
\definecolor{naranjauva}{RGB}{237,110,0}
\definecolor{verdecla}{rgb}{0.1,0.8,0.2} 
\definecolor{verdetlp}{RGB}{33,120,68}
\definecolor{granate}{rgb}{0.6,0,0.3} 
\definecolor{azulclarouva}{RGB}{130,115,186}
\definecolor{moradouva}{RGB}{186,31,181}
\definecolor{moradopod}{RGB}{97,45,98} 
\definecolor{moradomuyclaropod}{RGB}{160,129,161} 
\definecolor{moradoclaropod}{RGB}{149,78,153} 
\definecolor{verdeclaropod}{RGB}{151,194,184} 

\definecolor{verdeuva}{RGB}{122,154,1}
\definecolor{verdeosc}{RGB}{46, 139, 87} 
\definecolor{similar1}{RGB}{88,6,124}
\definecolor{similar2}{RGB}{186,14,0}
\definecolor{opuesto1}{RGB}{0,142,25}
\definecolor{opuesto2}{RGB}{178,185,0}
\definecolor{optotal}{RGB}{109,172,0}
\definecolor{amarillouva}{RGB}{172,104,43}     
\definecolor{rosapalo}{rgb}{0.83,0.48,0.38}   
\definecolor{violeta}{rgb}{0.7,0.5,0.8}       
\definecolor{rojo}{RGB}{219,0,0}                

\makeatletter
\DeclareRobustCommand\widecheck[1]{{\mathpalette\@widecheck{#1}}}
\def\@widecheck#1#2{%
    \setbox\z@\hbox{\m@th$#1#2$}%
    \setbox\tw@\hbox{\m@th$#1%
       \widehat{%
          \vrule\@width\z@\@height\ht\z@
          \vrule\@height\z@\@width\wd\z@}$}%
    \dp\tw@-\ht\z@
    \@tempdima\ht\z@ \advance\@tempdima2\ht\tw@ \divide\@tempdima\thr@@
    \setbox\tw@\hbox{%
       \raise\@tempdima\hbox{\scalebox{1}[-1]{\lower\@tempdima\box
\tw@}}}%
    {\ooalign{\box\tw@ \cr \box\z@}}}
\makeatother

\begin{document}

\title{Surjectivity of the asymptotic Borel map
in Carleman-Roumieu ultraholomorphic classes defined by regular sequences}
\author{Javier Jim\'enez-Garrido, Javier Sanz and Gerhard Schindl}
\date{\today}

\maketitle

{ \small
\begin{center}
{\bf Abstract}
\end{center}

We study the surjectivity of, and the existence of right inverses for, the asymptotic Borel map in Carleman-Roumieu ultraholomorphic classes defined by regular sequences in the sense of E.~M.~Dyn'kin. We extend previous results by J.~Schmets and M.~Valdivia, by V.~Thilliez, and by the authors, and show the prominent role played by an index associated with the sequence and introduced by Thilliez.
The techniques involve regular variation, integral transforms and characterization results of A. Debrouwere in a half-plane, steming from his study of the surjectivity of the moment mapping in general Gelfand-Shilov spaces.

}

\bigskip

\section{Introduction}

The concept of asymptotic expansion, introduced by H. Poincar\'e in 1886, has played an essential role in the understanding of the analytical meaning of the formal power series solutions to large classes of functional equations (ordinary and partial differential equations, difference and $q$-difference equations, and so on).
The existence of such an expansion for a complex holomorphic function in a sector $S$ of the Riemann surface of the logarithm amounts to a precise control on the growth of its derivatives, and this fact gives the link with ultraholomorphic classes, on whose elements' derivatives are usually imposed local or global bounds in terms of a weight sequence $\M=(M_n)_{n\in\N_0}$ of positive real numbers. See Subsection~\ref{subsectCarlemanclasses} for an account in this respect. The asymptotic Borel map sends a function in one of such classes into its formal power series of asymptotic expansion, and in many instances it is important to decide about its injectivity and surjectivity when considered between suitable spaces. We refer the reader to our previous paper~\cite{JimenezSanzSchindlInjectSurject}, whose introduction contains a non comprehensive historical account of the results in this respect, and where the problem of injectivity in unbounded sectors and for general weight sequences is completely closed, by solving a pending case not covered by the powerful results of S. Mandelbrojt~\cite{Mandelbrojt} and B. Rodr{\'\i}guez-Salinas~\cite{Salinas}.\par

Regarding surjectivity, the classical Borel-Ritt-Gevrey theorem of B.~Malgrange and J.-P.~Ramis~\cite{Ramis1}, solving the case of Gevrey asymptotics, was extended to different more general situations by J.~Schmets and M.~Valdivia~\cite{SchmetsValdivia00}, V.~Thilliez~\cite{Thilliez95,Thilliez03} and the authors~\cite{SanzFlatProxOrder,JimenezSanzSchindlInjectSurject}. For a weight sequence $\M$, our main satisfactory results have been the following:
\begin{itemize}
\item[(i)] The Borel map is never bijective~\cite[Theorem 3.17]{JimenezSanzSchindlInjectSurject}.
\item[(ii)] The strong nonquasianalyticity condition (equivalent to the fact that the index $\ga(\M)$ of Thilliez is positive) is necessary for surjectivity~\cite[Lemma 4.5]{JimenezSanzSchindlInjectSurject}.
\item[(iii)] For a sector $S_{\ga}$ of opening $\pi\gamma$ ($\ga>0$) and under the hypothesis of moderate growth for $\M$, surjectivity has been characterized (at least for uniform asymptotics, and except for the limiting case in some situations) by the condition $\ga<\ga(\M)$~\cite[Theorem 4.17]{JimenezSanzSchindlInjectSurject}.
\item[(iv)] Surjectivity was characterized whenever $\M$ admits a nonzero proximate order~\cite[Theorem 6.1]{SanzFlatProxOrder}.
\end{itemize}

The present paper intends to go one step further and complete the partial information given in~\cite[Theorem 4.14]{JimenezSanzSchindlInjectSurject} concerning the case of regular weight sequences in the sense of E.~M.~Dyn'kin~\cite{Dynkin80}, which instead of moderate growth satisfy the milder condition of derivation closedness (see Subsection~\ref{subsectstrregseq} for the precise definitions). Moreover, the existence of extension operators, right inverses for the Borel map, is studied in this general case. It is interesting to note that the condition $(\beta_2)$, introduced by H.-J.~Petzsche~\cite{Pet} in a similar study for ultradifferentiable classes, plays again a prominent role here, and its relationship with other conditions of rapid variation is elucidated. In particular, the condition $\ga(\M)=\infty$, stronger than $(\beta_2)$, guarantees the surjectivity of the Borel map and the existence of global extension operators for any sector in the Riemann surface of the logarithm.

We have not considered in this paper the closely related case of Beurling ultraholomorphic classes. The surjectivity of the Borel map in this setting, for $\ga<\ga(\M)$ and under the moderate growth condition, was established by V. Thilliez~\cite[Cor. 3.4.1]{Thilliez03}, and A. Debrouwere~\cite{DebrouwereExtOperBeurling} has very recently proved the existence of extension operators under the same hypotheses by using results from the splitting theory of Fréchet spaces.

\section{Preliminaries}\label{sectPrelimin}
\subsection{Notation}

We set $\N:=\{1,2,...\}$, $\N_{0}:=\N\cup\{0\}$.
$\mathcal{R}$ stands for the Riemann surface of the logarithm, where the notation $z=|z|e^{i\theta}$ refers to the element $(|z|,\theta)\in(0,\infty)\times\R$.
$\C[[z]]$ is the space of formal power series in $z$ with complex coefficients.

For $\gamma>0$, we consider unbounded sectors bisected by direction 0,
$$S_{\gamma}:=\{z\in\mathcal{R}:|\hbox{arg}(z)|<\frac{\gamma\,\pi}{2}\}$$
or, in general, bounded or unbounded sectors
$$S(d,\alpha,r):=\{z\in\mathcal{R}: |\hbox{arg}(z)-d|<\frac{\alpha\,\pi}{2},\ |z|<r\},\quad
S(d,\alpha):=\{z\in\mathcal{R}:|\hbox{arg}(z)-d|<\frac{\alpha\,\pi}{2}\}$$
with bisecting direction $d\in\R$, opening $\alpha\,\pi$ and (in the first case) radius $r\in(0,\infty)$.

A sector $T$ is said to be a \emph{proper subsector} of a sector $S$ if $\overline{T}\subset S$ (where the closure of $T$ is taken in $\mathcal{R}$, and so the vertex of the sector is not under consideration). In case such $T$ is also bounded, we say it is a \emph{bounded proper subsector} of $S$.

\subsection{Weight sequences and their properties}\label{subsectstrregseq}

In what follows, $\bM=(M_p)_{p\in\N_0}$ will always stand for a sequence of positive real numbers, and we will always assume that $M_0=1$. We define its {\it sequence of quotients}  $\m=(m_p)_{p\in\N_0}$ by
$m_p:=\frac{M_{p+1}}{M_p}$, $p\in \N_0$; clearly, the knowledge of $\M$ amounts to that of $\m$, since $M_p=m_0\cdots m_{p-1}$, $p\in\N$. We will denote by small letters the quotients of a sequence given by the corresponding capital letters. The following properties for a sequence will play a role in this paper:

\begin{itemize}
\item[(i)]  $\M$ is \emph{logarithmically convex} (for short, (lc)) if
$$M_{p}^{2}\le M_{p-1}M_{p+1},\qquad p\in\N.$$%
\item[(ii)] $\M$ is \emph{stable under differential operators} or satisfies the \emph{derivation closedness condition} (briefly, (dc)) if there exists $D>0$ such that $$M_{p+1}\leq D^{p+1} M_{p},  \qquad p\in\N_{0}. $$
\item[(iii)]  $\M$ is of, or has, \emph{moderate growth} (briefly, (mg)) whenever there exists $A>0$ such that
$$M_{p+q}\le A^{p+q}M_{p}M_{q},\qquad p,q\in\N_0.$$
\item[(iv)]  $\M$ satisfies the condition (snq) if there exists $B>0$ such that
$$
\sum^\infty_{q= p}\frac{M_{q}}{(q+1)M_{q+1}}\le B\frac{M_{p}}{M_{p+1}},\qquad p\in\N_0.$$
\end{itemize}

It will be convenient to introduce the notation $\hM:=(p!M_p)_{p\in\N_0}$.
All these properties are preserved when passing from $\M$ to $\hM$. In the classical work of H.~Komatsu~\cite{komatsu}, the properties (lc), (dc) and (mg) are denoted by $(M.1)$, $(M.2)'$ and $(M.2)$, respectively, while (snq) for $\M$ is the same as property $(M.3)$ for $\widehat{\M}$.
Obviously, (mg) implies (dc).

The sequence of quotients $\bm$ is nondecreasing if and only if $\M$ is (lc). In this case, it is well-known that $(M_p)^{1/p}\leq m_{p-1}$ for every $p\in\N$, the sequence $((M_p)^{1/p})_{p\in\N}$ is nondecreasing, and $\lim_{p\to\infty} (M_p)^{1/p}= \infty$ if and only if $\lim_{p\to\infty} m_p= \infty$. In order to avoid trivial situations, we will restrict from now on to
(lc) sequences $\M$ such that $\lim_{p\to\infty} m_p =\infty$, which will be called \emph{weight sequences}.
It is immediate that if $\hM$ satisfies (lc) and $\M$ satisfies (snq), then $\M$ is a weight sequence.

Following E.~M.~Dyn'kin~\cite{Dynkin80}, if $\M$ is a weight sequence and satisfies (dc), we say $\hM$ is \emph{regular}. According to V.~Thilliez~\cite{Thilliez03}, if $\M$ satisfies (lc), (mg) and (snq), we say $\M$ is \emph{strongly regular}; in this case $\M$ is a weight sequence, and the corresponding $\hM$ is regular.

We mention some interesting examples. In particular, those in (i) and (iii) appear in the applications of summability theory to the study of formal power series solutions for different kinds of equations.
\begin{itemize}
\item[(i)] The sequences $\M_{\a,\b}:=\big(p!^{\a}\prod_{m=0}^p\log^{\b}(e+m)\big)_{p\in\N_0}$, where $\a>0$ and $\b\in\R$, are strongly regular (in case $\b<0$, the first terms of the sequence have to be suitably modified in order to ensure (lc)). In case $\b=0$, we have the best known example of strongly regular sequence, $\M_{\a}:=\M_{\a,0}=(p!^{\a})_{p\in\N_{0}}$, called the \emph{Gevrey sequence of order $\a$}.

\item[(ii)] The sequence $\M_{0,\b}:=(\prod_{m=0}^p\log^{\b}(e+m))_{p\in\N_0}$, with $\b>0$, satisfies (lc) and (mg), and $\bm$ tends to infinity, but (snq) is not satisfied.
\item[(iii)] For $q>1$, $\M_q:=(q^{p^2})_{p\in\N_0}$ satisfies (lc), (dc) and (snq), but not (mg).
\end{itemize}

Two sequences $\bM=(M_{p})_{p\in\N_0}$ and $\bL=(L_{p})_{p\in\N_0}$ of positive real numbers, with respective quotients $\bm$ and $\bl$, are said to be:
\begin{itemize}
\item[(i)] \emph{equivalent}, and we write $\M\approx\bL$, if there exist positive constants $A,B$ such that
$$A^pM_p\le L_p\le B^pM_p,\qquad p\in\N_0.$$
\item[(ii)] \emph{strongly equivalent}, and we write $\bm\simeq\bl$, if there exist positive constants $a,b$ such that
$$am_p\le \ell_p\le bm_p,\qquad p\in\N_0.$$
\end{itemize}
Whenever $\bm\simeq\bl$ we have $\M\approx\bL$, but not conversely.

As an example, for $\a>0$ we set $\bL_{\a}:=(\Ga(1+\a p))_{p\in\N_0}$, where $\Ga$ denotes the Eulerian Gamma function; it is well-known that $\bl_{\a}\simeq((p+1)^{\a})_{p\in\N_0}$ and so $\bL_{\a}\approx\M_{\a}$, the Gevrey sequence of order $\a$.

Conditions (dc) and (mg) are clearly preserved by $\approx$, and so also by $\simeq$, for general sequences; (snq) is obviously preserved for weight sequences by $\simeq$, but also by $\approx$ (see the work of H.-J. Petzsche~\cite[Cor. 3.2]{Pet} for an indirect argument, and our paper~\cite[Cor. 3.14]{JimenezSanzSchindlIndex} for a direct proof of a more general statement).

Given two sequences $\M$ and $\bL$, we use the notation $\M\cdot\bL=(M_nL_n)_{n\in\N_0}$ and $\M/\bL=(M_n/L_n)_{n\in\N_0}$.
We will use the fact that $\M$ satisfies (mg), respectively (dc), if and only if $\M\cdot\bL_{\a}$ or $\M/\bL_{\a}$ satisfy (mg), resp. (dc), for some $\a>0$.

\subsection{Asymptotic expansions, ultraholomorphic classes and the asymptotic Borel map}\label{subsectCarlemanclasses}

In this paragraph $S$ is a sector and $\M$ a  sequence. We start recalling the concept of asymptotic expansion.

We say a holomorphic function $f$ in $S$ admits the formal power series $\widehat{f}=\sum_{p=0}^{\infty}a_{p}z^{p}\in\C[[z]]$ as its $\{\M\}$-\emph{asymptotic expansion} in $S$ (when the variable tends to 0) if for every bounded proper subsector $T$ of $S$ there exist $C_T,A_T>0$ such that for every $p\in\N_0$, one has
\begin{equation*}\Big|f(z)-\sum_{n=0}^{p-1}a_nz^n \Big|\le C_TA_T^pM_{p}|z|^p,\qquad z\in T.
\end{equation*}
If the expansion exists, it is unique, and we will write $f\sim_{\{\bM\}}\widehat{f}$ in $S$. $\widetilde{\mathcal{A}}_{\{\M\}}(S)$ stands for the space of functions admitting $\{\M\}$-asymptotic expansion in $S$.\par

We say a holomorphic function $f:S\to\C$ admits $\widehat{f}$ as its \emph{uniform $\{\M\}$-asymptotic expansion in $G$ (of type $1/A$ for some $A>0$)} if there exists $C>0$ such that for every $p\in\N_0$, one has
\begin{equation}\Big|f(z)-\sum_{n=0}^{p-1}a_nz^n \Big|\le CA^pM_{p}|z|^p,\qquad z\in S.\label{desarasintunifo}
\end{equation}
In this case we write $f\sim_{\{\M\},A}^u\widehat{f}$ in $S$, and $\widetilde{\mathcal{A}}^u_{\{\M\},A}(S)$ denotes the space of functions admitting uniform $\{\M\}$-asymptotic expansion of type $1/A$ in $S$, endowed with the norm
$$
\left\|f\right\|_{\M,A,\overset{\sim}{u}}:=\sup_{z\in S,n\in\N_{0}}\frac{|f(z)-\sum_{k=0}^{n-1}a_kz^k|}{A^{n}M_{n}|z|^n},
$$
which makes it a Banach space. $\widetilde{\mathcal{A}}^u_{\{\M\}}(S)$ stands for the $(LB)$ space of functions admitting uniform $\{\M\}$-asymptotic expansion in $S$, obtained as the union of the previous classes when $A$ runs over $(0,\infty)$.
When the type needs not be specified, we simply write $f\sim_{\{\M\}}^u\widehat{f}$ in $S$.
Note that, taking $p=0$ in~\eqref{desarasintunifo}, we deduce that every function in $\widetilde{\mathcal{A}}^u_{\{\M\}}(S)$ is a bounded function.

Finally, we define for every $A>0$ the class $\mathcal{A}_{\{\M\},A}(S)$ consisting of the functions holomorphic in $S$ such that
$$
\left\|f\right\|_{\M,A}:=\sup_{z\in S,n\in\N_{0}}\frac{|f^{(n)}(z)|}{A^{n}M_{n}}<\infty.
$$
($\mathcal{A}_{\{\M\},A}(S),\left\|\,\cdot\, \right\| _{\M,A}$) is a Banach space, and $\mathcal{A}_{\{\M\}}(S):=\cup_{A>0}\mathcal{A}_{\{\M\},A}(S)$ is called a \emph{Carleman-Roumieu ultraholomorphic class} in the sector $S$, whose natural inductive topology makes it an $(LB)$ space.\par
We warn the reader that these notations do not agree with the ones used in~\cite{SanzFlatProxOrder,JimenezSanzSchindlInjectSurject}, where $\widetilde{\mathcal{A}}_{\{\M\}}(S)$ was denoted by $\widetilde{\mathcal{A}}_{\M}(S)$,
$\widetilde{\mathcal{A}}^u_{\{\M\}}(S)$ by
$\widetilde{\mathcal{A}}^u_{\M}(S)$, $\mathcal{A}_{\{\M\},A}(S)$ by $\mathcal{A}_{\M/\bL_1,A}(S)$, and $\mathcal{A}_{\{\M\}}(S)$ by $\mathcal{A}_{\M/\bL_1}(S)$.

If $\M$ is (lc), the spaces $\mathcal{A}_{\{\M\}}(S)$, $\widetilde{\mathcal{A}}^u_{\{\M\}}(S)$ and $\widetilde{\mathcal{A}}_{\{\M\}}(S)$ are algebras, and if $\M$ is (dc) they are stable under taking derivatives.
Moreover, if $\M\approx\bL$ the corresponding classes coincide.

Since the derivatives of $f\in\mathcal{A}_{\{\M\},A}(S)$ are Lipschitzian, for every $n\in\N_{0}$ one may define
\begin{equation}\label{eq.deriv.at.0.def}
f^{(p)}(0):=\lim_{z\in S,z\to0 }f^{(p)}(z)\in\C.
\end{equation}

As a consequence of Taylor's formula and Cauchy's integral formula for the derivatives, there is a close relation between Carleman-Roumieu ultraholomorphic classes and the concept of asymptotic expansion (the proof may be easily adapted from~\cite{balserutx}).

\begin{prop}\label{propcotaderidesaasin}
Let $\M$ be a sequence and $S$ be a sector. Then,
\begin{enumerate}[(i)]
\item If $f\in\mathcal{A}_{\{\hM\},A}(S)$ then $f$ admits $\widehat{f}:=\sum_{p\in\N_0}\frac{1}{p!}f^{(p)}(0)z^p$ as its uniform $\{\M\}$-asymptotic expansion in $S$ of type $1/A$, where $(f^{(p)}(0))_{p\in\N_0}$ is given by \eqref{eq.deriv.at.0.def}. Moreover, $\|f\|_{\M,A,\overset{\sim}{u}}\le \|f\|_{\hM,A}$, and so the identity $\mathcal{A}_{\{\hM\},A}(S)\hookrightarrow \widetilde{\mathcal{A}}^u_{\{\M\},A}(S)$ is continuous. Consequently, we also have that
$$\mathcal{A}_{\{\hM\}}(S)\subseteq \widetilde{\mathcal{A}}^u_{\{\M\}}(S) \subseteq  \widetilde{\mathcal{A}}_{\{\M\}}(S),$$
and $\mathcal{A}_{\{\hM\}}(S)\hookrightarrow \widetilde{\mathcal{A}}^u_{\{\M\}}(S)$ is continuous.

\item $f\in\widetilde{\mathcal{A}}_{\{\M\}}(S)$ if and only if for every (bounded or, if possible, unbounded) proper subsector $T$ of $S$ there exists $A_T>0$ such that $f|_T\in \mathcal{A}_{\{\hM\},A_T}(T)$.
In case any of the previous holds and $f\sim_{\{\M\}}\sum^\infty_{p=0} a_p z^p$, then for every such $T$ and every $p\in\N_0$ one has
\begin{equation}\label{equaCoeffAsympExpanLimitDeriv}
a_p=\lim_{ \genfrac{}{}{0pt}{}{z\to0}{z\in T}} \frac{f^{(p)}(z)}{p!},
\end{equation}
and we can set ${f^{(p)}(0)}:=p!a_p$.

\item If $S$ is unbounded and $T$ is a proper (bounded or unbounded) subsector of $S$, then there exists a constant $c=c(T,S)>0$ such that the restriction to $T$, $f|_T$, of functions $f$ defined on $S$ and admitting uniform $\{\M\}$-asymptotic expansion in $S$ of type $1/A>0$, belongs to $\mathcal{A}_{\{\hM\},cA}(T)$, and $\|f|_T\|_{\hM,cA}\le \|f\|_{\M,A,\overset{\sim}{u}}$. So, the restriction map from $\widetilde{\mathcal{A}}^u_{\{\M\},A}(S)$ to $\mathcal{A}_{\{\hM\},cA}(T)$ is continuous, and it is also continuous from $\widetilde{\mathcal{A}}^u_{\{\M\}}(S)$ to $\mathcal{A}_{\{\hM\}}(T)$.
\end{enumerate}
\end{prop}

One may accordingly define classes of formal power series
$$\C[[z]]_{\{\M\},A}=\Big\{\widehat{f}=\sum_{n=0}^\infty a_nz^n\in\C[[z]]:\, \left|\,\widehat{f} \,\right|_{\M,A}:=\sup_{p\in\N_{0}}\displaystyle \frac{|a_{p}|}{A^{p}M_{p}}<\infty\Big\}.$$%
$(\C[[z]]_{\{\M\},A},\left| \,\cdot\,  \right|_{\M,A})$ is a Banach space and we put $\C[[z]]_{\{\M\}}:=\cup_{A>0}\C[[z]]_{\{\M\},A}$, again an $(LB)$ space.

Given $f\in\widetilde{\mathcal{A}}_{\{\M\}}(S)$ with $f\sim_{\{\M\}}\widehat{f}$, and taking into account~\eqref{equaCoeffAsympExpanLimitDeriv}, it is straightforward that $\widehat{f}\in\C[[z]]_{\{\M\}}$, so it is natural to consider the \textit{asymptotic Borel map}
$$
\widetilde{\mathcal{B}}:\widetilde{\mathcal{A}}_{\{\M\}}(S)\longrightarrow \C[[z]]_{\{\M\}}$$
sending a function $f\in\widetilde{\mathcal{A}}_{\{\M\}}(S)$ into its $\{\M\}$-asymptotic expansion $\widehat{f}$.
By Proposition~\ref{propcotaderidesaasin}.(i)
the asymptotic Borel map may be defined in $\widetilde{\mathcal{A}}^u_{\{\M\}}(S)$, $\mathcal{A}_{\{\hM\}}(S)$ and $\mathcal{A}_{\{\hM\},A}(S)$ (in the last case, with target space $\C[[z]]_{\{\M\},A}$).

We would like to highlight that, alternatively, the target space for the Borel map could be considered to be a space of sequences comprising the derivatives at 0 of a function $f$ in the classes, as defined in~\eqref{eq.deriv.at.0.def}, and subject to the corresponding control on the growth of their terms. This equivalent approach has been followed by many authors, and in particular in the works of J.~Schmets and M.~Valdivia~\cite{SchmetsValdivia00} and A.~Debrouwere~\cite{momentsdebrouwere}. Note that their results, stated in this paper as Theorems~\ref{teor.Andreas},~\ref{th.ExtensOperHalfplane} and~\ref{th.SchmetsValdiviafastgrowth}, have been adapted to our setting.

If $\M$ is (lc), $\widetilde{\mathcal{B}}$ is a homomorphism of algebras; if $\M$ is also (dc), differentiation commutes with $\widetilde{\mathcal{B}}$. Moreover, it is continuous when considered between the corresponding Banach or $(LB)$ spaces previously introduced. Finally, note that if $\M\approx\bL$, then $\C[[z]]_{\{\M\}}=\C[[z]]_{\{\bL\}}$, and the corresponding Borel maps are in all cases identical.

Since the problem under study is invariant under rotation, we will focus on the surjectivity of the Borel map in unbounded sectors $S_{\ga}$. So, we define
\begin{align*}
S_{\{\hM\}}:=&\{\ga>0; \quad \widetilde{\mathcal{B}}:\mathcal{A}_{\{\hM\}}(S_\ga)\longrightarrow \C[[z]]_{\{\M\}} \text{ is surjective}\} ,\\
\widetilde{S}^u_{\{\M\}}:=&\{\ga>0; \quad\widetilde{\mathcal{B}}:\widetilde{\mathcal{A}}^u_{\{\M\}}(S_\ga)\longrightarrow \C[[z]]_{\{\M\}} \text{ is surjective}\}, \\
\widetilde{S}_{\{\M\}}:=&\{\ga>0;\quad \widetilde{\mathcal{B}}:\widetilde{\mathcal{A}}_{\{\M\}}(S_\ga)\longrightarrow \C[[z]]_{\{\M\}} \text{ is surjective} \}.
\end{align*}
We again note that these intervals were respectively denoted by $S_{\M}$, $\widetilde{S}^u_{\M}$ and $\widetilde{S}_{\M}$ in~\cite{JimenezSanzSchindlInjectSurject}.\par

It is clear that $S_{\{\hM\}}$, $\widetilde{S}^u_{\{\M\}}$ and $\widetilde{S}_{\{\M\}}$ are either empty or left-open intervals  having $0$ as endpoint, called \emph{surjectivity intervals}.
Using~Proposition~\ref{propcotaderidesaasin}, items (i) and (iii),
we easily see that
\begin{align}
(\widetilde{S}^u_{\{\M\}})^{\circ}\subseteq S_{\{\hM\}}\subseteq \widetilde{S}^u_{\{\M\}} \subseteq\widetilde{S}_{\{\M\}},
\label{equaContentionSurjectIntervals}
\end{align}
where $I^{\circ}$ stands for the interior of the interval $I$.

\section{Surjectivity results for regular sequences}\label{sectSurjectivity}

In the study of the surjectivity the index $\ga(\M)$, introduced in this regard by V.~Thilliez~\cite[Sect.\ 1.3]{Thilliez03} for strongly regular sequences $\M$, will play a central role.
His definition makes sense for (lc) sequences, in this case $\ga(\M)\in[0,\infty]$, and it may be equivalently expressed by different conditions:
\begin{enumerate}[(i)]
\item A sequence $(c_p)_{p\in\N_0}$ is \emph{almost increasing} if there exists $a>0$ such that for every $p\in\N_0$ we have that $c_p\leq a c_q $ for every $ q\geq p$.
It was proved in~\cite{JimenezSanzSRSPO,JimenezSanzSchindlIndex} that for any weight sequence $\M$ one has
\begin{equation}\label{equa.indice.gammaM.casicrec}
\ga(\bM)=\sup\{\ga>0:(m_{p}/(p+1)^\ga)_{p\in\N_0}\hbox{ is almost increasing} \}.
\end{equation}
\item For any $\b>0$ we say that $\m$ satisfies the condition $(\ga_{\b})$ if there exists $A>0$ such that
\begin{equation}\label{equa.condgammabeta}
\sum^\infty_{\ell=p} \frac{1}{(m_\ell)^{1/\b}}\leq \frac{A (p+1) }{(m_p)^{1/\b}},  \qquad p\in\N_0.\tag{$\ga_\b$}
\end{equation}
Using this condition, which was introduced for $\beta=1$ by H. Komatsu~\cite{komatsu} (and named $(\ga_{1})$ after H.-J. Petzsche~\cite{Pet}), and generalized for $\b\in\N$ by J. Schmets and M. Valdivia~\cite{SchmetsValdivia00}, we can obtain (see~\cite{PhDJimenez,JimenezSanzSchindlIndex}) that
\begin{equation}\label{equa.indice.gammaM.gamma_r}
\ga(\M)=\sup\{\b>0; \,\, \m \,\, \text{satisfies} \,\, \eqref{equa.condgammabeta}\,\};\quad \ga(\M)>\b\iff\m\text{ satisfies }(\ga_\b).
\end{equation}
\end{enumerate}

Whenever $\widehat{\M}=(p!M_p)_{p\in\N_0}$ is (lc) we have (see~\cite[Ch.~2]{PhDJimenez} and \cite[Cor.~3.13]{JimenezSanzSchindlIndex}) that
$\ga(\M)>0$ if and only if $\M$ is (snq).
We recall also the following result for later use.

\begin{lemma}[\cite{JimenezSanzSchindlIndex}, Remark 3.15]\label{lemma.gammaMgreaterthan1}
For an arbitrary sequence $\M_1$ such that $\gamma(\M_1)>1$,
there exists a weight sequence $\M_2$ such that  $\widehat{\bm}_2\simeq\bm_1$, and so $\ga(\widehat{\M}_2)=\ga(\M_1)$.
\end{lemma}

A straightforward verification shows that for any sequence $\M$ and for every $s>0$ one has
\begin{align}
\ga((p!^sM_p)_{p\in\N_0})&=\ga((\Ga(1+sp)M_p)_{p\in\N_0})=\ga(\M)+s,\label{equa.gamma.producto}\\
\ga((M_p/p!^s)_{p\in\N_0})&=\ga(M_p/(\Ga(1+sp))_{p\in\N_0})=\ga(\M)-s.\label{equa.gamma.cociente}
\end{align}

As a consequence of the characterization of the surjectivity of the Borel map in the ultradifferentiable setting given by H.-J. Petzsche\cite[Thm. 3.5]{Pet}, we proved the following result, already announced by V. Thilliez in~\cite{Thilliez03}.

\begin{lemma}[\cite{JimenezSanzSchindlInjectSurject}, Lemma 4.5]\label{lemmaSurjectivityImpliessnq}
Let $\M$ be a weight sequence. If $\widetilde{S}_{\{\M\}}\neq\emptyset$, then $\M$ has (snq) or, equivalently, $\ga(\M)>0$.
\end{lemma}

Our aim in this section is to solve (except for some limiting cases) the problem of surjectivity whenever $\M$ is a weight sequence satisfying (dc) or, in other words, $\hM$ is a regular sequence in the sense of Dyn'kin. Our previous main result is the following. We denote by $\lfloor x\rfloor$ the greatest integer not exceeding $x$.

\begin{theo}[\cite{JimenezSanzSchindlInjectSurject}, Thm. 4.14 and Cor. 4.15]\label{th.SurjectivityUniformAsymp.plus.dc}
Let $\M$ be a weight sequence satisfying (dc).
\begin{itemize}
\item[(i)] Let $\alpha>0$ be such that $\widetilde{\mathcal{B}}:\widetilde{\mathcal{A}}^u_{\{\M\}}(S_{\alpha})\to \C[[z]]_{\{\M\}}$ is surjective. Then, $\ga(\M)>\lfloor \alpha \rfloor $.
\item[(ii)] If we have that $\widetilde{S}^u_{\{\M\}}=(0,\infty)$, then $S_{\{\hM\}} =\widetilde{S}^u_{\{\M\}} =\widetilde{S}_{\{\M\}}=(0,\infty)$ and $\ga(\M)=\infty$.
\end{itemize}
One has $S_{\{\hM\}}\subseteq \widetilde{S}^u_{\{\M\}} \subseteq (0,\lfloor \ga(\M) \rfloor +1)$; if moreover $\ga(\M)\in\N$, then $S_{\{\hM\}}\subseteq \widetilde{S}^u_{\{\M\}} \subseteq (0,\ga(\M))$.
\end{theo}

At that moment and to the best of our knowledge, no general surjectivity result had been proved for regular $\hM$, except for the special case of the $q$-Gevrey sequences $\M_q=(q^{p^2})_{p\in\N_0}$, $q>1$, see C.~Zhang~\cite{Zhang}. In a recent collaboration of the first two authors with A. Debrouwere~\cite{DebrouwereJimenezSanzStieltjesmoment} we have studied the existence and uniqueness of solutions for the Stieltjes moment  problem in Gelfand-Shilov spaces, subspaces of the Schwartz space of rapidly decreasing smooth functions for which the growth of the products of monomials times the derivatives of their elements is controlled in terms of weight sequences. By a suitable application of the Fourier transform, there exists a close connection between this problem and the surjectivity or injectivity of the asymptotic Borel map in ultraholomorphic classes in a half-plane, and so our results in~\cite{JimenezSanzSchindlInjectSurject} could be transferred, providing a complete solution for the surjectivity of the moment map whenever strongly regular sequences are considered, and only a partial one for regular sequences. The key point for our coming results is a new work by A. Debrouwere~\cite{momentsdebrouwere}, where the surjectivity of the Stieltjes moment problem for regular sequences has been characterized by using only functional-analytic methods. Again thanks to the Fourier transform (but in the opposite direction) he has taken this information into the asymptotic framework. We state next a version adapted to our needs: firstly, while we ask for $\M$ to be (lc), it is enough that $\hM$ is; secondly, the condition $\ga(\M)>1$ amounts, in view of~\eqref{equa.gamma.producto} and~\eqref{equa.indice.gammaM.gamma_r}, to the fact that $\hM$ satisfies $(\ga_2)$, which is the condition appearing in~\cite[Thm. 7.4.(b)]{momentsdebrouwere}.

\begin{theo}[\cite{momentsdebrouwere}]\label{teor.Andreas}
Let $\hM$ be regular. The following are equivalent:
\begin{itemize}
\item[(i)] $\widetilde{\mathcal{B}}\colon \mathcal{A}_{\{\hM\}}(S_1)\to\C[[z]]_{\{\M\}}$ is surjective.
\item[(ii)] $\ga(\M)>1$.
\end{itemize}
\end{theo}

We highlight that (i)$\Rightarrow$(ii) is slightly weaker than part (i) of Theorem~\ref{th.SurjectivityUniformAsymp.plus.dc} when $\a=1$; on the other hand, the implication (ii)$\Rightarrow$(i) provides the first general surjectivity result for weight sequences not subject to condition (mg) (apart from a result of J.~Schmets and M.~Valdivia for rapidly varying sequences which we will comment on later).

However, the previous method seems to be valid only for a half-plane. We will be able to carry the information to the case of a general sector by applying general Laplace, $\mathcal{L}_{\a}$, and Borel, $\mathcal{B}_{\a}$, transforms of order $\a>0$, which basically arise from the classical transforms (inverse of each other) combined with ramifications of exponent $\a$.
Namely, we will follow the approach in Sections 5.5 and 5.6 of the book of W. Balser~\cite{balserutx}. We recall that, for $0<\a<2$, one considers the Laplace kernel function
$$
e_{\a}(z):=\frac{1}{\a}z^{1/\a}\exp(-z^{1/\a}),\qquad z\in S_{\a},
$$
whose moment function is
$$
m_{\a}(\lambda):=\int_{0}^{\infty}t^{\lambda-1}e_{\a}(t)dt=\Gamma(1+\a \lambda), \qquad\Re(\lambda)\ge 0,
$$
and the corresponding Borel kernel function
$$E_{\a}(z):=\sum_{n=0}^\infty \frac{z^n}{m_{\a}(n)}=
\sum_{n=0}^\infty \frac{z^n}{\Gamma(1+\a n)},\qquad z\in\C,$$
which is the classical Mittag-Leffler function of order $\a$.

Subsequently, given a function $f$ holomorphic in a sector $S=S(d,\alpha)$ and with suitable growth, for any direction $\tau$ in $S$ the \textit{$\a$-Laplace transform in direction $\tau$} of $f$ is defined as
\begin{equation*}
(\mathcal{L}_{\a,\tau}f)(z):=\int_0^{\infty(\tau)}e_{\a}(u/z)f(u)\frac{du}{u},\quad
|\arg(z)-\tau|<\a\pi/2,\ |z|\textrm{ small enough},
\end{equation*}
where the integral is taken along the half-line parameterized by $t\in(0,\infty)\mapsto te^{i\tau}$. The family $\{\mathcal{L}_{\a,\tau}f\}_{\tau\textrm{\,in\,}S}$ defines a holomorphic function $\mathcal{L}_{\a}f$ named the \textit{$\a$-Laplace transform} of $f$.

Secondly, let $S=S(d,\b,r)$ be a sector with $\b>\a$, and $f:S\to \C$ be holomorphic in $S$ and continuous at 0 (i.e. the limit of $f$ at 0 exists when $z$ tends to 0 in every proper subsector of $S$).
For $\tau\in\R$ such that $|\tau-d|<(\b-\a)\pi/2$ we may consider a path
$\delta_{\a}(\tau)$ in $S$ like the ones used in the classical Borel transform, consisting of a segment from the origin to a point $z_0$ with $\arg(z_0)=\tau+\a(\pi+\varepsilon)/2$ (for some
suitably small $\varepsilon\in(0,\pi)$), then the circular arc $|z|=|z_0|$ from $z_0$ to
the point $z_1$ on the ray $\arg(z)=\tau-\a(\pi+\varepsilon)/2$ (traversed clockwise), and
finally the segment from $z_1$ to the origin.

The \textit{$\a$-Borel transform in direction $\tau$} of $f$ is then defined as
$$
(\mathcal{B}_{\a,\tau}f)(u):=\frac{-1}{2\pi i}\int_{\delta_{\a}(\tau)}E_{\a}(u/z)f(z)\frac{dz}{z},\quad
u\in S(\tau,\varepsilon_0), \quad \varepsilon_0\textrm{ small enough}.
$$
The family
$\{\mathcal{B}_{\a,\tau}f\}_{\tau}$ defines the \textit{$\a$-Borel transform} of $f$, holomorphic in the sector $S(d,\b-\a)$ and denoted by $\mathcal{B}_{\a}f$.\par

The formal $\a$-Laplace and $\a$-Borel transforms, defined from $\C[[z]]$ into $\C[[z]]$, are respectively given by
$$\widehat{\mathcal{L}}_{\a}\big(\sum_{p=0}^{\infty}a_{p}z^{p}\big):= \sum_{p=0}^{\infty}\Gamma(1+\a p)a_{p}z^{p},\qquad \widehat{\mathcal{B}}_{\a}\big(\sum_{p=0}^{\infty}a_{p}z^{p}\big):= \sum_{p=0}^{\infty}\frac{a_{p}}{\Gamma(1+\a p)}z^{p}.$$

The following result, involving two sequences, can be found in a slightly different form in~\cite[Thms. 27 and 28]{balserutx}, where only the case of two Gevrey sequences is considered, and in~\cite[Thm. 3.16]{lastramaleksanzJMAA15}, where a general sequence and a sequence admitting a nonzero proximate order intervene. Here, we consider an intermediate situation.

\begin{theo}\label{teorrelacdesartransfBL}
Suppose $\M$ is an arbitrary sequence, and $\a,\ga>0$. Let $f\in\widetilde{\mathcal{A}}^{u}_{\{\bM\}}(S_\ga)$
and $f\sim^u_{\{\M\}}\widehat{f}$. Recall that
$\bL_{\a}:=(\Ga(1+\a p))_{p\in\N_0}$. Then, the following hold:
\begin{itemize}
\item[(i)] For every $\b$ with $0<\b<\ga$ one has $$\mathcal{L}_{\a}f\in \widetilde{\mathcal{A}}^{u}_{\{\bM\cdot\bL_{\a}\}}(S_{\b+\a}) \quad\textrm{and}\quad  \mathcal{L}_{\a}f\sim^u_{\{\M\cdot\bL_{\a}\}} \widehat{\mathcal{L}}_{\a}\widehat{f}.
    $$
    Moreover, there exist $C,c>0$, depending only on $\a$, $\b$ and $\ga$, such that for every $A>0$ and every $f\in\widetilde{\mathcal{A}}^{u}_{\{\bM\},A}(S_\ga)$ one has
    $\|\mathcal{L}_{\a}f\|_{\M\cdot\bL_{\a},cA,\overset{\sim}{u}}\le C\|f\|_{\M,A,\overset{\sim}{u}}$, and so the maps $\mathcal{L}_{\a}\colon \widetilde{\mathcal{A}}^{u}_{\{\bM\},A}(S_\ga)\to \widetilde{\mathcal{A}}^{u}_{\{\bM\cdot\bL_{\a}\},cA}(S_{\b+\a})$ and $\mathcal{L}_{\a}\colon \widetilde{\mathcal{A}}^{u}_{\{\bM\}}(S_\ga)\to \widetilde{\mathcal{A}}^{u}_{\{\bM\cdot\bL_{\a}\}}(S_{\b+\a})$ are continuous.
\item[(ii)] Suppose $\ga>\a$. For every $\b$ with $\a<\b<\ga$ one has
    $$\mathcal{B}_{\a}f\in \widetilde{\mathcal{A}}^{u}_{\{\bM/\bL_{\a}\}}(S_{\b-\a})\quad \textrm{and}\quad \mathcal{B}_{\a}f\sim^u_{\{\M/\bL_{\a}\}} \widehat{\mathcal{B}}_{\a}\widehat{f}.
    $$
    Moreover, there exist $C,c>0$, depending only on $\a$, $\b$ and $\ga$, such that for every $A>0$ and every $f\in\widetilde{\mathcal{A}}^{u}_{\{\bM\},A}(S_\ga)$ one has
    $\|\mathcal{B}_{\a}f\|_{\M/\bL_{\a},cA,\overset{\sim}{u}}\le C\|f\|_{\M,A,\overset{\sim}{u}}$, and so the maps $\mathcal{B}_{\a}\colon \widetilde{\mathcal{A}}^{u}_{\{\bM\},A}(S_\ga)\to \widetilde{\mathcal{A}}^{u}_{\{\bM/\bL_{\a}\},cA}(S_{\b-\a})$ and $\mathcal{B}_{\a}\colon \widetilde{\mathcal{A}}^{u}_{\{\bM\}}(S_\ga)\to \widetilde{\mathcal{A}}^{u}_{\{\bM/\bL_{\a}\}}(S_{\b-\a})$ are continuous.
\end{itemize}
\end{theo}

With the help of this result we can complete the information in~\eqref{equaContentionSurjectIntervals}. We will use the clear fact that the formal Laplace and Borel transforms, $\widehat{\mathcal{L}}_{\a}$ and $\widehat{\mathcal{B}}_{\a}$, are (topological) isomorphisms between the space $\C[[z]]_{\{\M\}}$ and $\C[[z]]_{\{\M\cdot\bL_{\a}\}}$, respectively $\C[[z]]_{\{\M/\bL_{\a}\}}$, for an arbitrary sequence $\M$.

\begin{lemma}\label{lemma.secondintercontclosurethird}
For any weight sequence $\M$, $\widetilde{S}_{\{\M\}}$ is contained in the closure of $\widetilde{S}^u_{\{\M\}}$ in $(0,\infty)$.
\end{lemma}
\begin{proof1}
Due to the form of these intervals, it is equivalent to prove that whenever $\ga>0$ belongs to $\widetilde{S}_{\{\M\}}$, one has $(0,\ga)\subseteq \widetilde{S}^u_{\{\M\}}$. Let us see that any $\b\in(0,\ga)$ belongs to $\widetilde{S}^u_{\{\M\}}$. Choose positive real numbers $\a,\b'$ such that $\a<\b<\b'<\ga$. First, we deduce that $\widetilde{\mathcal{B}}\colon \widetilde{\mathcal{A}}^u_{\{\M/\bL_{\a}\}}(S_{\b'-\a}) \to\C[[z]]_{\{\M/\bL_{\a}\}}$ is surjective. Given $\widehat{g}\in\C[[z]]_{\{\M/\bL_{\a}\}}$, we know $\widehat{f}:=\widehat{\mathcal{L}}_{\a}\widehat{g}\in\C[[z]]_{\{\M\}}$.
Since $\widetilde{\mathcal{B}}\colon \widetilde{\mathcal{A}}_{\{\M\}}(S_\ga)\to\C[[z]]_{\{\M\}}$ is surjective, there exists $f\in\widetilde{\mathcal{A}}_{\{\M\}}(S_\ga)$ such that $f\sim_{\{\M\}}\widehat{f}$. One may apply the Borel transform $\mathcal{B}_{\a}$ to $f$, and the proof of~\cite[Thm. 3.16.(ii)]{lastramaleksanzJMAA15} shows that from the asymptotic estimates in bounded proper subsectors of $S_{\ga}$ for $f$ one can deduce uniform asymptotic estimates in $S_{\b'-\a}$ for $\mathcal{B}_{\a}f$, and moreover $\mathcal{B}_{\a}f\sim^{u}_{\{\M/\bL_{\a}\}}\widehat{g}$, as desired.

Subsequently, a similar use of the Laplace transform $\mathcal{L}_{\a}$ shows, by taking into account Theorem~\ref{teorrelacdesartransfBL}.(i), that $\widetilde{\mathcal{B}}\colon \widetilde{\mathcal{A}}^u_{\{\M\}}(S_{\b})\to\C[[z]]_{\{\M\}}$ is also surjective, and we conclude.
\end{proof1}

We can now state our first main result.

\begin{theo}\label{teorSurject.dc}
Let $\hM$ be a regular sequence such that $\ga(\M)>0$. Then,
$$(0,\ga(\M))\subseteq S_{\{\hM\}}\subseteq \widetilde{S}^u_{\{\M\}}\subseteq \widetilde{S}_{\{\M\}}\subseteq(0,\ga(\M)].
$$
\end{theo}
\begin{proof1}
According to~\eqref{equaContentionSurjectIntervals} and Lemma~\ref{lemma.secondintercontclosurethird}, it suffices to prove that $(0,\ga(\M))\subseteq\widetilde{S}^u_{\{\M\}}\subseteq(0,\ga(\M)]$.

Firstly, we suppose $0<\ga<\ga(\M)$ and prove that $\ga\in\widetilde{S}^u_{\{\M\}}$. We distinguish two cases:
\begin{itemize}
\item[(a.1)] If $\ga(\M)>1$, it suffices to work with $\ga>1$. Take $\ga'$ such that $\ga<\ga'<\ga(\M)$. The sequence $\mathbb{P}_1:=\hM/\bL_{\ga'-1}$ satisfies (dc) and, thanks to~\eqref{equa.gamma.cociente}, $\ga(\mathbb{P}_1)=\ga(\M)-\ga'+2>2$.
By Lemma~\ref{lemma.gammaMgreaterthan1}, there exists a weight sequence $\mathbb{P}_2$ such that $\widehat{\mathbb{P}_2}\approx\mathbb{P}_1$, $\ga(\mathbb{P}_2)=\ga(\mathbb{P}_1)-1>1$, and which satisfies (dc). Theorem~\ref{teor.Andreas} applies, so $\widetilde{\mathcal{B}}\colon \mathcal{A}_{\{\widehat{\mathbb{P}_2}\}}(S_1)\to\C[[z]]_{\{\mathbb{P}_2\}}$ is surjective, and the same holds when the map departs from $\widetilde{\mathcal{A}}^u_{\{\mathbb{P}_2\}}(S_1)$. Combining this fact with an application of the Laplace transform $\mathcal{L}_{\ga'-1}\colon \widetilde{\mathcal{A}}^{u}_{\{\mathbb{P}_2\}}(S_1)\to \widetilde{\mathcal{A}}^{u}_{\{\mathbb{P}_2\cdot\bL_{\ga'-1}\}}(S_{\ga})$, Theorem~\ref{teorrelacdesartransfBL}.(i) shows that, since $\ga<\ga'=1+(\ga'-1)$, also $\widetilde{\mathcal{B}}\colon \widetilde{\mathcal{A}}^u_{\{\mathbb{P}_2\cdot\bL_{\ga'-1}\}}(S_{\ga}) \to\C[[z]]_{\{\mathbb{P}_2\cdot\bL_{\ga'-1}\}}$ is surjective. We conclude by observing that $\mathbb{P}_2\cdot\bL_{\ga'-1}\approx(\mathbb{P}_1/\bL_1)\cdot\bL_{\ga'-1}=
\hM/\bL_1=\M$,
so that the corresponding classes coincide and $\ga\in\widetilde{S}^u_{\{\M\}}$.

\item[(a.2)] If $\ga(\M)\le 1$, choose $\a\in(0,1)$ such that $\ga(\M)+\a>1$. Hence, $\M\cdot\bL_{\a}$ is a weight sequence satisfying (dc) and, by using~\eqref{equa.gamma.producto}, $\ga(\M\cdot\bL_{\a})>1$. Given $\ga'$ such that $\ga<\ga'<\ga(\M)$, by the previous item (a.1) we know that $\ga'+\a\in\widetilde{S}^u_{\{\M\cdot\bL_{\a}\}}$. In this case, we may combine this fact with an application of the Borel transform $\mathcal{B}_{\a}\colon \widetilde{\mathcal{A}}^{u}_{\{\bM\cdot\bL_{\a}\}}(S_{\ga'+\a})\to \widetilde{\mathcal{A}}^{u}_{\{\bM\}}(S_{\ga})$, and Theorem~\ref{teorrelacdesartransfBL}.(ii) implies that, since $\ga<\ga'$, also $\ga\in\widetilde{S}^u_{\{\M\}}$, as desired.
\end{itemize}
Secondly, we take $\ga\in\widetilde{S}^u_{\{\M\}}$ and we will prove that $\ga\le\ga(\M)$. We again have different cases:
\begin{itemize}
\item[(b.1)] If $0<\ga<1$, consider positive real numbers $\a,\ga'$ with $1-\a<\ga'<\ga$. By applying the Laplace transform $\mathcal{L}_{\a}\colon \widetilde{\mathcal{A}}^{u}_{\{\M\}}(S_{\ga})\to \widetilde{\mathcal{A}}^{u}_{\{\M\cdot\bL_{\a}\}}(S_{\ga'+\a})$, Theorem~\ref{teorrelacdesartransfBL}.(i) shows that $\ga'+\a\in\widetilde{S}^u_{\{\M\cdot\bL_{\a}\}}$. Observe that $\ga'+\a>1$, so we deduce by restriction to the half-plane $S_1$ that, according to Proposition~\ref{propcotaderidesaasin}.(iii), also $1\in S_{\{\hM\cdot\bL_{\a}\}}$. Theorem~\ref{teor.Andreas} implies then that $\ga(\M\cdot\bL_{\a})>1$ or, equivalently by~\eqref{equa.gamma.producto}, $\ga(\M)>1-\a$. Since $\a$ can be chosen arbitrarily while keeping $1-\a<\ga$, we deduce $\ga(\M)\ge\ga$.
\item[(b.2)] If $\ga\in\N$, we know that $\ga(\M)>\ga$ by Theorem~\ref{th.SurjectivityUniformAsymp.plus.dc}.(i).
\item[(b.3)] If $\ga\in(1,\infty)\setminus\N$, again by Theorem~\ref{th.SurjectivityUniformAsymp.plus.dc}.(i) we deduce that $\ga(\M)>\lfloor\ga\rfloor$, so that the sequence $\mathbb{P}_1:=\hM/\bL_{\lfloor\ga\rfloor}$ is such that $\ga(\mathbb{P}_1)>1$ by using~\eqref{equa.gamma.producto} and~\eqref{equa.gamma.cociente}. Hence, from Lemma~\ref{lemma.gammaMgreaterthan1} there exists a weight sequence $\mathbb{P}_2$ such that $\mathbb{P}_2\approx\M/\bL_{\lfloor\ga\rfloor}$, and  $\mathbb{P}_2$ will also satisfy (dc). Consider a value $\ga'$ with $\lfloor\ga\rfloor<\ga'<\ga$. An application of the Borel transform $\mathcal{B}_{\lfloor\ga\rfloor}\colon \widetilde{\mathcal{A}}^{u}_{\{\bM\}}(S_{\ga})\to \widetilde{\mathcal{A}}^{u}_{\{\bM/\bL_{\lfloor\ga\rfloor}\}} (S_{\ga'-\lfloor\ga\rfloor})$ and Theorem~\ref{teorrelacdesartransfBL}.(ii) shows that $\ga'-\lfloor\ga\rfloor\in \widetilde{S}^u_{\{\M/\bL_{\lfloor\ga\rfloor}\}}$ or, equivalently, $\ga'-\lfloor\ga\rfloor\in \widetilde{S}^u_{\{\mathbb{P}_2\}}$. Since $ \ga'-\lfloor\ga\rfloor\in(0,1)$, we may invoke item (b.1) and deduce that $\ga(\mathbb{P}_2)\ge\ga'-\lfloor\ga\rfloor$, what amounts to $\ga(\M)\ge\ga'$. We conclude by making $\ga'$ tend to $\ga$.
\end{itemize}
\end{proof1}

The previous result confirms that, as indicated by V. Thilliez in~\cite[Sect. 3.3]{Thilliez03}, the moderate growth condition (mg) was of a technical nature for surjectivity. For weight sequences satisfying (dc) it is only pending to determine whether $\ga(\M)$ belongs or not to the surjectivity intervals. In the particular case that $\ga(\M)\in\N$ we know $S_{\{\hM\}} =\widetilde{S}^u_{\{\M\}}=(0,\ga(\M))$ (see Theorem~\ref{th.SurjectivityUniformAsymp.plus.dc}). On the other hand, if $\M$ admits a nonzero proximate order (which is more restrictive than strong regularity, but a common situation in applications) we also know that $\widetilde{S}_{\{\M\}}=(0,\ga(\M)]$, see~\cite[Thm. 6.1]{SanzFlatProxOrder} and~\cite[Thm. 4.24]{JimenezSanzSchindlInjectSurject}.

In general, for an arbitrary weight sequence we have no proof of surjectivity for any opening, and the surjectivity intervals could possibly be empty; however, in~\cite[Thm. 4.10 and Cor. 4.11]{JimenezSanzSchindlInjectSurject} we have obtained that $\widetilde{S}_{\{M\}}\subset(0,\lfloor\ga(\M)\rfloor +1]$.

In view of the previous information, our conjecture is that $S_{\{\hM\}} =\widetilde{S}^u_{\{\M\}}=(0,\ga(\M))$ and $\widetilde{S}_{\{\M\}}=(0,\ga(\M)]$ in general.

\section{Global extension operators}

One may ask about the existence of extension operators, right inverses for the asymptotic Borel map. This can be done, in principle, in the Banach spaces $\widetilde{\mathcal{A}}^u_{\{\M\},A}(S)$ and  $\mathcal{A}_{\{\M\},A}(S)$, which we call the local case, or in the $(LB)$ spaces  $\widetilde{\mathcal{A}}^u_{\{\M\}}(S)$ and $\mathcal{A}_{\{\M\}}(S)$, which we refer to as the global one. The first situation was studied by V. Thilliez, see~\cite[Thm. 3.2.1]{Thilliez03}, who obtained local extension operators with an scaling of the type for strongly regular sequences in sectors $S_{\ga}$ as long as $\ga<\ga(\M)$.

In the global situation and in the ultradifferentiable setting, H.-J.~Petzsche introduced the condition
\begin{equation}\label{equa.beta2}
\forall \varepsilon>0,\ \exists k\in\N,\ k>1:
\limsup_{p\to\infty} \left(\frac{M_{kp}}{M_p}\right)^{\frac{1}{(k-1)p}}\frac{1}{m_{kp-1}} \le\varepsilon,\tag{$\b_2$}
\end{equation}
which again appeared in the results of J.~Schmets and M.~Valdivia~\cite{SchmetsValdivia00} and A.~Debrouwere~\cite{momentsdebrouwere} about the existence of global extension operators in the ultraholomorphic framework. Please note that the sequence of quotients considered in these two previously cited papers results from our sequence $\bm$ after an index shift by 1, what explains the slightly different expression given here to condition $(\b_2)$.
We subsequently mention a version of the result by A.~Debrouwere adapted to our needs, in a similar way as in Theorem~\ref{teor.Andreas}.

\begin{theo}[\cite{momentsdebrouwere}, Thm. 7.4]\label{th.ExtensOperHalfplane}
Suppose $\hM$ is a regular sequence. The following are equivalent:
\begin{itemize}
\item[(i)] There exists a global extension operator $U_{\M}:\C[[z]]_{\{\M\}}\to\mathcal{A}_{\{\hM\}}(S_{1})$.
\item[(ii)] $\ga(\M)>1$, and $\M$ satisfies $(\b_2)$.
\end{itemize}
\end{theo}

The use of Laplace and Borel transforms of arbitrary positive order allows us to generalize this statement. We will also take into account  that condition $(\b_2)$ is evidently stable under strong equivalence $\simeq$ and, as a consequence of Stirling's formula (see~\cite[Lemma 2.2.(b)]{SchmetsValdivia00}), a sequence $\M$ satisfies $(\b_2)$ if and only if $\M\cdot\bL_{\a}$ or $\M/\bL_{\a}$ satisfies $(\b_2)$ for some/any $\a>0$.

\begin{theo}\label{th.GlobalExtOperRegSeq}
Suppose $\hM$ is a regular sequence, and let $r>0$. Each of the following statements implies the next one:
\begin{itemize}
\item[(i)] $r<\ga(\M)$, and $\M$ satisfies $(\b_2)$.
\item[(ii)] There exists a global extension operator $U_{\M,r}:\C[[z]]_{\{\M\}}\to\mathcal{A}_{\{\hM\}}(S_{r})$.
\item[(iii)] There exists a global extension operator $V_{\M,r}:\C[[z]]_{\{\M\}}\to\widetilde{\mathcal{A}}^u_{\{\M\}}(S_{r})$.
\item[(iv)] $r\le \ga(\M)$, and $\M$ satisfies $(\b_2)$.
\end{itemize}
\end{theo}

\begin{proof1}
(i)$\implies$(ii) We consider two cases:
\begin{itemize}
\item[(a.1)] Suppose $r>1$, and take a real number $r'$ with $r<r'<\ga(\M)$. Reasoning as in the proof of Theorem~\ref{teorSurject.dc}.(a.1), there exists a weight sequence $\mathbb{P}$
    such that $\boldsymbol{p}\simeq\bm/\bl_{r'-1}$, satisfies (dc) and $(\b_2)$, and $\ga(\mathbb{P})=\ga(\M)+1-r'>1$. Theorem~\ref{th.ExtensOperHalfplane} provides an extension operator $U\colon\C[[z]]_{\{\M/\bL_{r'-1}\}}\to \mathcal{A}_{\{\hM/\bL_{r'-1}\}}(S_{1})$. By Proposition~\ref{propcotaderidesaasin}.(i), this induces an extension operator $\widetilde{U}\colon\C[[z]]_{\{\M/\bL_{r'-1}\}}\to \widetilde{\mathcal{A}}^u_{\{\M/\bL_{r'-1}\}}(S_{1})$. Theorem~\ref{teorrelacdesartransfBL}.(i) implies that the composition $\mathcal{L}_{r'-1}\circ\widetilde{U}\circ \widehat{\mathcal{B}}_{r'-1}$ will be an extension operator from $\C[[z]]_{\{\M\}}$ to $\widetilde{\mathcal{A}}^u_{\{\M\}}(S_{\rho})$ for every $0<\rho<r'=1+(r'-1)$. If we choose $\rho=(r+r')/2>r$, the restriction of the elements of this last space to $S_r$ provides, by Proposition~\ref{propcotaderidesaasin}.(iii), the extension operator $U_{\M,r}:\C[[z]]_{\{\M\}}\to\mathcal{A}_{\{\hM\}}(S_{r})$ we were looking for.
\item[(a.2)] If $r\le 1$, consider $\a$ such that $\a+r>1$, and take $r'$ with $r<r'<\ga(\M)$. The sequence $\M\cdot\bL_{\a}$ satisfies $(\b_2)$ and $\ga(\M\cdot\bL_{\a})>r'+\a>1$. By item (a.1), there exists an extension operator $U\colon\C[[z]]_{\{\M\cdot\bL_{\a}\}}\to \mathcal{A}_{\{\hM\cdot\bL_{\a}\}}(S_{r'+\a})$. Again Proposition~\ref{propcotaderidesaasin}.(i) allows us to obtain an extension operator $\widetilde{U}\colon\C[[z]]_{\{\M\cdot\bL_{\a}\}}\to \widetilde{\mathcal{A}}^u_{\{\M\cdot\bL_{\a}\}}(S_{r'+\a})$. Now, Theorem~\ref{teorrelacdesartransfBL}.(ii) implies that  $\mathcal{B}_{\a}\circ\widetilde{U}\circ \widehat{\mathcal{L}}_{\a}$ will be an extension operator from $\C[[z]]_{\{\M\}}$ to, say, $\widetilde{\mathcal{A}}^u_{\{\M\}}(S_{(r+r')/2})$, and the restriction of the elements of this space to $S_{r}$ provides the desired extension operator as before.
\end{itemize}\par\noindent
(ii)$\implies$(iii) Obvious from Proposition~\ref{propcotaderidesaasin}.(i).\par\noindent
(iii)$\implies$(iv) We consider again two cases:
\begin{itemize}
\item[(b.1)] Suppose $r>1$, and take a real number $r'$ with $1<r'<r$. The existence of $V_{\M,r}$ implies that the corresponding Borel map is surjective in $S_r$, and by Theorem~\ref{teorSurject.dc} we have $\ga(\M)\ge r$. So, repeating the argument in (a.1), there exists a weight sequence $\mathbb{P}$ such that $\boldsymbol{p}\simeq\bm/\bl_{r'-1}$, satisfies (dc) and $\ga(\mathbb{P})=\ga(\M)+1-r'>1$. Since the classes associated with $\M$ and $\mathbb{P}\cdot\bL_{r'-1}$ agree, we have an extension operator $\widetilde{V}_{\M,r}:\C[[z]]_{\{\mathbb{P}\cdot\bL_{r'-1}\}} \to\widetilde{\mathcal{A}}^u_{\{\mathbb{P}\cdot\bL_{r'-1}\}}(S_{r})$. Note that $1+(r-r')/2<r-(r'-1)$, and so the mapping $\mathcal{B}_{r'-1}\circ \widetilde{V}_{\M,r}\circ \widehat{\mathcal{L}}_{r'-1}$ is an extension operator from
    $\C[[z]]_{\{\mathbb{P}\}}$ to $\widetilde{\mathcal{A}}^u_{\{\mathbb{P}\}}(S_{1+(r-r')/2})$. The restriction of the elements of this last space to $S_1$ provides, by Proposition~\ref{propcotaderidesaasin}.(iii), an extension operator $U_{\M,r}:\C[[z]]_{\{\mathbb{P}\}}\to \mathcal{A}_{\{\widehat{\mathbb{P}}\}}(S_{1})$.
    Then, Theorem~\ref{th.ExtensOperHalfplane} guarantees that $\mathbb{P}$ satisfies $(\b_2)$, and so $\M$ will also do according to the stability properties of $(\b_2)$. Moreover, $\ga(\mathbb{P})>1$, from where $\ga(\M)>r'$. Since $r'$ was arbitrarily close to $r$, we deduce that $\ga(\M)\ge r$, as desired.
\item[(b.2)] If $r\le 1$, consider $\a$ such that $\a+r>1$, and take $\a'>\a$. Since $r+\a<r+\a'$, Theorem~\ref{teorrelacdesartransfBL}.(i) asserts that the mapping $\mathcal{L}_{\a'}\circ V_{\M,r}\circ \widehat{\mathcal{B}}_{\a'}$ will be an extension operator from $\C[[z]]_{\{\M\cdot\bL_{\a'}\}}$ to $\widetilde{\mathcal{A}}^u_{\{\M\cdot\bL_{\a'}\}}(S_{r+\a})$. We can apply item (b.1) and deduce that $\M\cdot\bL_{\a'}$ satisfies $(\b_2)$, and so $\M$ will also do, and that  $\ga(\M\cdot\bL_{\a'})=\ga(\M)+\a'\ge r+\a$. We conclude by making $\a'$ tend to $\a$.
\end{itemize}
\end{proof1}

Our conjecture is that (i), (ii) and (iii) in Theorem~\ref{th.GlobalExtOperRegSeq} are equivalent, but we are not able to fill the gap at this moment.

Observe that if $\M$ is a weight sequence satisfying $(\b_2)$, we may apply Lemma 2.4 in~\cite{SchmetsValdivia00} to the sequence $\hM$ and deduce that $\ga(\M)>0$. So, if $\hM$ is regular and satisfies $(\b_2)$, one can always obtain extension operators for $0<r<\ga(\M)$ thanks to the previous theorem.

In the last part of our study, we want to determine the weight sequences for which extension operators exist for sectors of arbitrary opening. In this respect, J. Schmets and V. Valdivia state the following result for sequences with fast growth. Please recall that the sequence of quotients considered by these authors results from our sequence $\bm$ after an index shift by 1.

\begin{theo}[\cite{SchmetsValdivia00}, Thm. 5.6]\label{th.SchmetsValdiviafastgrowth}
Let $\M$ be a weight sequence such that
\begin{equation}\label{equa.cond.SchmetsValdivia}
\text{for every }r\in\N,\ (m_{n-1}/n^r)_{n\in\N}\text{ is increasing from some term on.}
\end{equation}
The following are equivalent:
\begin{itemize}
\item[(i)] For every $r\in\N$, there exists a global extension operator $U_{\M,r}:\C[[z]]_{\{\M\}}\to\mathcal{A}_{\{\hM\}}(S_{r})$.
\item[(ii)] For some $r\in\N$, there exists a global extension operator $U_{\M,r}:\C[[z]]_{\{\M\}}\to\mathcal{A}_{\{\hM\}}(S_{r})$.
\item[(iii)] $\M$ satisfies $(\b_2)$.
\end{itemize}
\end{theo}

However, it turns out that the conditions~\eqref{equa.cond.SchmetsValdivia} and $(\b_2)$ are related to each other. The connection among these and other conditions of fast growth, usually appearing in the literature, can be inferred from the theory of rapid variation (see the classical book of Bingham et al.~\cite{BingGoldTeug89}) and our study of the indices and orders of regular variation associated with weight sequences~\cite{JimenezSanzSchindlIndex}.

We recall that in the study of the injectivity of the Borel map for ultraholomorphic classes in unbounded sectors, completed in~\cite{JimenezSanzSchindlInjectSurject}, the growth index (introduced in~\cite{SanzFlatProxOrder}, see also~\cite{JimenezSanzSRSPO})
\begin{equation}\label{equa.defi.indexomegaM}
\o(\M):= \displaystyle\liminf_{n\to\infty} \frac{\log(m_{n})}{\log(n)}\in[0,\infty]
\end{equation}
played a prominent role. Moreover, the moderate growth condition (mg) is satisfied by $\M$ precisely when the upper Matuszewska index associated with its sequence of quotients, $\a(\m)$, is finite (see~\cite[Cor. 3.17]{JimenezSanzSchindlIndex}), and we recall that, for a general weight sequence,
\begin{equation}\label{equa.Ineq.Indices}
0\le\ga(\M)=\b(\m)\le \o(\M)\le\a(\m)\le\infty
\end{equation}
always holds, where $\b(\m)$ is the lower Matuszewska index associated with $\m$ (\cite[Rem. 3.4 and Thm. 3.10]{JimenezSanzSchindlIndex}).

\begin{prop}\label{propRelationshipConditRapidVariation}
Let $\M$ be a weight sequence. Each of the following statements implies the next one, and only the implications (ii)$\implies$(iii)$\implies$(iv) may be reversed:
\begin{itemize}
\item[(i)] $\M$ satisfies~\eqref{equa.cond.SchmetsValdivia}.
\item[(ii)] $\ga(\M)=\infty$.
\item[(iii)] For every $k\in\N$, $k\ge 2$, one has $\lim_{n\to\infty}\displaystyle\frac{m_{kn}}{m_n}=\infty$.
\item[(iv)] There exists $k_0\in\N$, $k_0\ge 2$, such that $\lim_{n\to\infty}\displaystyle\frac{m_{k_0n}}{m_n}=\infty$.
\item[(v)] $\M$ satisfies $(\beta_2)$.
\item[(vi)] $\lim_{n\to\infty}\displaystyle\frac{m_{n}}{M_n^{1/n}}=\infty$.
\item[(vii)] $\o(\M)=\infty$.
\item[(viii)] $\a(\m)=\infty$ (in other words, $\M$ does not satisfy (mg)).
\end{itemize}
\end{prop}
\begin{proof1}
\noindent (i)$\implies$(ii) The condition~\eqref{equa.cond.SchmetsValdivia} clearly implies that the sequence $(m_{n-1}/n^r)_{n\in\N}$ is almost increasing for every $r\in\N$. As indicated in~\cite[Rem.~3.8]{JimenezSanzSchindlIndex}, this entails the same for the sequence $(m_{n}/n^r)_{n\in\N}$, and we only need to recall~\eqref{equa.indice.gammaM.casicrec} in order to deduce $\ga(\M)=\infty$. On the contrary, consider the sequence $\M$ whose quotients $(m_n)_{n\in\N_0}$ are given by
$$
m_{n-1}=\begin{cases}
q^{2n+1} &\text{ if }n\neq 2^k+1\text{ for every }k\in\N_0,\\
q^{2n-1} &\text{ if }n=2^k+1\text{ for some }k\in\N_0,
\end{cases}
$$
where $q>1$ and $n\ge 1$. It is not difficult to check that the sequence $(m_{n-1}/n)_{n\in\N}$ is not eventually increasing, while $(m_n/n^r)_{n\in\N}$ is almost increasing for every $r\in\N$, and so $\ga(\M)=\infty$.\par\noindent
(ii)$\implies$(iii) Consider the function $f(x):=m_{\lfloor x\rfloor}$, $x\ge 1$, which is measurable, nondecreasing, and whose lower Matuszewska index $\b(f)$ equals that of $\m$, which is precisely $\ga(\M)=\infty$ (see~\cite{BingGoldTeug89}, \cite[Sect. 3]{JimenezSanzSchindlIndex}). This means that $f$ belongs to the class $MR_{\infty}$ of rapid variation (\cite[p. 83]{BingGoldTeug89}), what, by Proposition 2.4.4.(iii) in~\cite{BingGoldTeug89}, amounts to the fact that $\lim_{n\to\infty}m_{\lfloor \lambda x\rfloor}/m_{\lfloor x\rfloor}=\infty$ for every $\lambda>1$. This implies (iii).\par\noindent
(iii)$\implies$(iv) Obvious.\par\noindent
(iv)$\implies$(ii) If (iv) is satisfied, for $k\in\N$ with $k\ge k_0$ and for every $\b>0$ we have $$\liminf_{p\to\infty}m_{kp}/(k^{\b}m_p)=\infty,
$$
and so also $\lim_{k\to\infty}\liminf_{p\to\infty}m_{kp}/(k^{\b}m_p)=\infty$. By Theorem 3.11 in~\cite{JimenezSanzSchindlIndex} we see that $\ga(\M)>\b$ and, $\b$ being arbitrary, we deduce that $\ga(\M)=\infty$.\par\noindent
(iv)$\implies$(v) We put $\widetilde{m}_n:=m_{n-1}$ for $n\in\N$, $\widetilde{m}_0:=1$, and $\widetilde{\m}:=(\widetilde{m}_n)_{n\in\N_0}$. By the previous items, (iv) amounts to $\ga(\M)=\infty=\b(\m)$, and this is the same as $\b(\widetilde{\m})=\infty$ again by~\cite[Rem.~3.8]{JimenezSanzSchindlIndex}. This, in turn, is equivalent to the fact that condition (iv) is satisfied now by the sequence $\widetilde{\m}$, and this appears in the work of H.-J. Petzsche~\cite{Pet} as condition $(\b_2^0)$. Proposition 1.6.(a) in~\cite{Pet} proves that $(\b_2^0)$ for $\widetilde{\m}$ implies (v).

The implication cannot be reversed because of the Example 1.8.(a) in~\cite{Pet}, and the fact that $\m$ satisfies (iv) if and only if $\widetilde{\m}$ does.\par\noindent
(v)$\implies$(vi) In~\cite[p. 304]{Pet} it is proved that (v) implies that $\lim_{n\to\infty}\displaystyle\frac{m_{n-1}}{M_n^{1/n}}=\infty$, and (vi) follows since $m_n\ge m_{n-1}$, $n\in\N$.

On the contrary, by Example 1.8.(b) in~\cite{Pet} there exists an increasing sequence of positive real numbers $\widetilde{\m}:=(\widetilde{m}_n)_{n\in\N_0}$ with $\widetilde{m}_0=1$ such that, if one defines $m_n=\widetilde{m}_{n+1}$, $n\in\N_0$, and the corresponding $\M$, it holds that $\M$ does not satisfy $(\b_2)$ and
$\lim_{n\to\infty}\displaystyle\frac{m_{n-1}}{M_n^{1/n}}=\infty$, so that (vi) is satisfied.\par\noindent
(vi)$\implies$(vii) For convenience, we put $\a_n:=\log(m_n)$, $n\in\N_0$; $\b_0:=\a_0$, $\b_n:=\log\left(\frac{m_n}{M^{1/n}_n}\right)$, $n\geq1$.
By Lemma 3.8 in~\cite{JimenezSanzSRSPO} we know that
$$
\a_n=\sum^{n-1}_{k=0}\frac{\b_k}{k+1}+\b_n, \qquad n\in\N_0,
$$
and so, using (vi), we have
\begin{equation*}
\lim_{n\to\infty} \frac{(\a_{n+1}-\b_{n+1})-(\a_n-\b_n)}{\log(n+1)-\log(n)}=\lim_{n\to\infty} \frac{\b_n/(n+1)}{1/n}=\infty.
\end{equation*}
We deduce by Stolz's criterion that
\begin{equation*}
\lim_{n\to\infty} \frac{\a_n-\b_n}{\log(n)}=\infty,
\end{equation*}
and since $\b_n\ge 0$ for every $n$, we obtain that $\lim_{n\to\infty} \a_n/\log(n)=\infty$. The conclusion follows from~\eqref{equa.defi.indexomegaM}.

On the contrary, consider two increasing sequences of natural numbers $(q_k)_{k\in\N}$ and $(p_k)_{k\in\N}$ such that $q_1=2$, for every $k\in\N$ one puts $p_k:=q_k^2-1$, and the $q_k$ are recursively chosen under the condition
\begin{equation}\label{equa.conddefisubseqexample}
q_{k+1}\ge k\log(q_{k+1})+p_k,\quad k\ge 1.
\end{equation}
Note that one has $q_k<p_k<q_{k+1}$ for every $k\in\N$.
We define the sequence $\M$ whose quotients are given by
\begin{equation*}
m_{p-1}=\exp\left(\sum_{j=1}^p\delta_j\right),\;\;\;p\ge 1,
\end{equation*}
where $(\delta_j)_{j\ge 1}$ is the sequence with
$\delta_1=\delta_2=0$ and
$$\delta_j=\begin{cases}
0,&\text{ if }q_k+1\le j\le p_k\text{ for some $k\ge 1$,}\\
1,&\text{ if }p_k+1\le j\le q_{k+1}\text{ for some $k\ge 1$.}
\end{cases}$$
Since the $\delta_j$ are nonnegative, $\M$ automatically satisfies (lc). It is plain to show that for every $p\in\N$ one has
\begin{equation*}
\frac{m_{p}}{(M_p)^{1/p}}= \exp\left(\frac{1}{p}\sum_{j=1}^pj\delta_{j+1}\right).
\end{equation*}
Then,
\begin{align*}
\log\left(\frac{m_{p_k}}{(M_{p_k})^{1/p_k}}\right)&
=\frac{1}{p_k}\sum_{j=1}^{p_k}j\delta_{j+1}\le \frac{1}{p_k}\sum_{j=1}^{q_k}j=\frac{q_k}{q_k-1}\frac{1}{2}\le 1,
\end{align*}
and thus (vi) is violated.
On the other hand, note that for every $k\ge 1$ one has
$$\frac{\log(m_{q_{k+1}-1})}{\log(q_{k+1})}= \frac{1}{\log(q_{k+1})}\sum_{j=1}^{q_{k+1}}\delta_j\ge \frac{q_{k+1}-p_k}{\log(q_{k+1})}\ge k,$$
where the last inequality stems from~\eqref{equa.conddefisubseqexample}. From here,
$$\frac{\log(m_{p_{k+1}})}{\log(p_{k+1}+1)}= \frac{1}{\log(p_{k+1}+1)}\sum_{j=1}^{p_{k+1}+1}\delta_j= \frac{1}{\log(q_{k+1}^2)}\sum_{j=1}^{q_{k+1}}\delta_j\ge\frac{k}{2}.$$
A simple study of the monotonicity of the sequence $(\frac{\log(m_{p})}{\log(p+1)})_{p\in\N}$ implies that
$\lim_{p\rightarrow\infty}\frac{\log(m_{p})}{\log(p+1)}=+\infty$ (in particular, $\M$ is a weight sequence), and so $\omega(\M)=\infty$.
\par\noindent

(vii)$\implies$(viii) The implication comes from~\eqref{equa.Ineq.Indices}. However, from the theory of rapid variation we learn that strict inequalities are possible in every case in~\eqref{equa.Ineq.Indices}. A particular example showing that $\a(\m)=\infty$ and $\o(\M)<\infty$ may simultaneously hold can be found in~\cite[p. 106]{PhDJimenez}, resting on another example by M. Langenbruch~\cite{Langenbruch}.
\end{proof1}

As a first consequence, note that for strongly regular sequences surjectivity does hold for small openings and local extension operators exist with an scaling in the type (see~\cite[Thm.~3.2.1]{Thilliez03}), but no global extension operator is possible, since condition $(\b_2)$ avoids moderate growth.\par

Secondly, the next result clarifies the situation for rapidly growing sequences and avoids to impose the condition (dc).
Note that $\ga(\M)=\infty$ guarantees that (snq) is satisfied, but is independent from condition (dc).

\begin{theo}
Let $\M$ be a weight sequence. The following are equivalent:
\begin{itemize}
\item[(i)] $\ga(\M)=\infty$.
\item[(ii)] For every $r>0$, there exists a global extension operator $U_{\M,r}:\C[[z]]_{\{\M\}}\to\mathcal{A}_{\{\hM\}}(S_{r})$.
\item[(iii)] For every $r>0$, there exists a global extension operator $V_{\M,r}:\C[[z]]_{\{\M\}}\to\widetilde{\mathcal{A}}^u_{\{\M\}}(S_{r})$.
\item[(iv)] All the surjectivity intervals are $(0,\infty)$.
\end{itemize}
\end{theo}

\begin{proof1}
\noindent (i)$\implies$(ii) Given $r>0$, consider $r_0:=\lfloor r\rfloor +1>r$. From~\eqref{equa.gamma.producto} it is clear that $\ga(\hM)=\infty$, and by Proposition~\ref{propRelationshipConditRapidVariation} we have that $\hM$ satisfies $(\b_2)$. Moreover, \eqref{equa.indice.gammaM.gamma_r}  implies that $\widehat{\bm}$ satisfies $(\ga_{r_0+1})$. We can apply Theorem 5.4 in~\cite{SchmetsValdivia00} for the sequence $\hM$ and the positive integer $r_0+1$, and subsequently Theorem 5.5 in~\cite{SchmetsValdivia00} for the value $\a=r$, in order to obtain an extension operator $U_{\M,r}:\C[[z]]_{\{\M\}}\to\mathcal{A}_{\{\hM\}}(S_{r})$.\par\noindent
(ii)$\implies$(iii) It is clear by Proposition~\ref{propcotaderidesaasin}.(i).\par\noindent
(iii)$\implies$(iv) By the definition of global extension operators as right inverses for the Borel map, we obviously have $\widetilde{S}^u_{\{\M\}}=(0,\infty)$. Then, \eqref{equaContentionSurjectIntervals} leads to the statement.\par\noindent
(iv)$\implies$(i) It suffices to apply Theorem 4.10 in~\cite{JimenezSanzSchindlInjectSurject}.

\end{proof1}

We note that pathological situations are possible. For example, if $\hM$ is regular, $\ga(\M)\in(0,\infty)$ and~\eqref{equa.beta2} holds, we have surjectivity in $S_\ga$, with global right inverses, for every $\ga<\ga(\M)$, but surjectivity fails for $\ga>\ga(\M)$; since~\eqref{equa.beta2} implies $\o(\M)=\infty$, injectivity will not hold in any (narrow or wide) sector.

\vskip.2cm
\noindent\textbf{Acknowledgements}: The first two authors are partially supported by the Spanish Ministry of Economy, Industry and Competitiveness under the project MTM2016-77642-C2-1-P. The third author is supported by FWF-Project J~3948-N35, as a part of which he has been an external researcher at the Universidad de Valladolid (Spain) for the period October 2016 - December 2018, and by FWF-Project P32905.\par

\vskip.5cm

\noindent\textbf{Affiliation}:\\
J.~Jim\'{e}nez-Garrido:\\
Departamento de Matem\'aticas, Estad{\'\i}stica y Computaci\'on\\
Universidad de Cantabria\\
Facultad de Ciencias, Avda. de los Castros, 48, 39005 Santander, Spain.\\
Instituto de Investigaci\'on en Matem\'aticas IMUVA, Universidad de Va\-lla\-do\-lid\\
E-mail: jesusjavier.jimenez@unican.es.
\\
\vskip.5cm
\noindent J.~Sanz:\\
Departamento de \'Algebra, An\'alisis Matem\'atico, Geometr{\'\i}a y Topolog{\'\i}a\\
Universidad de Va\-lla\-do\-lid\\
Facultad de Ciencias, Paseo de Bel\'en 7, 47011 Valladolid, Spain.\\
Instituto de Investigaci\'on en Matem\'aticas IMUVA, Universidad de Va\-lla\-do\-lid\\
E-mail: jsanzg@am.uva.es.
\\
\vskip.5cm
\noindent G.~Schindl:\\
Fakult\"at f\"ur Mathematik, Universit\"at Wien,
Oskar-Morgenstern-Platz~1, A-1090 Wien, Austria.\\
E-mail: gerhard.schindl@univie.ac.at.
\end{document}